# The Laguerre process and generalized Hartman–Watson law

NIZAR DEMNI

*Laboratoire de Probabilités et Modèles Aléatoires, Université de Paris VI, 4 Place Jussieu, Case 188, F-75252 Paris Cedex 05. E-mail: demni@ccr.jussieu.fr*

In this paper, we study complex Wishart processes or the so-called Laguerre processes $(X_t)_{t \geq 0}$. We are interested in the behaviour of the eigenvalue process; we derive some useful stochastic differential equations and compute both the infinitesimal generator and the semi-group. We also give absolute-continuity relations between different indices. Finally, we compute the density function of the so-called generalized Hartman–Watson law as well as the law of $T_0 := \inf\{t, \det(X_t) = 0\}$ when the size of the matrix is 2.

*Keywords:* generalized Hartman–Watson law; Gross–Richards formula; Laguerre process; special functions of matrix argument

## 1. Introduction

The Real Wishart process is a symmetric matrix-valued process which was introduced by Bru [2] as follows. Let $B_t = (B_{ij}(t))_{i,j}$ be an $n \times m$ Brownian matrix and define $X_t = B_t^{\mathrm{T}} B_t$. The process $(X_t)_{t \geq 0}$ satisfies the following stochastic differential equation (SDE):

$$\begin{aligned} \mathrm{d}X_t &= B_t^{\mathrm{T}} \, \mathrm{d}B_t + \mathrm{d}B_t^{\mathrm{T}} B_t + n I_m \, \mathrm{d}t \\ &= \sqrt{X_t} \mathrm{d}N_t + \mathrm{d}N_t^{\mathrm{T}} \sqrt{X_t} + n I_m \, \mathrm{d}t, \qquad X_0 = B_0^{\mathrm{T}} B_0, \end{aligned}$$

where $I_m$ denotes the unit matrix, the superscript $^{\mathrm{T}}$ stands for the transpose, $\sqrt{X_t}$ is the matrix square root of the positive matrix $X_t$ and $(N_t)_{t \geq 0}$ is an $m \times m$ Brownian matrix. This process is called the Wishart process of dimension $n$, of size $m$, starting from $X_0$, and is denoted by $W(n, m, X_0)$. Then $W(\delta, m, X_0)$, where $\delta$ runs over the Gindikin ensemble $\{1, \ldots, m-1\} \cup \,]m-1, \infty[$, is defined as the unique solution of the latter SDE with $\delta$ instead of $n$. Thus, it can be viewed as an extension of the squared Bessel process to higher dimensions. In this way, Donati-Martin *et al.* [6] tried to extend some well-known properties of the squared Bessel processes to the matrix case and derived expressions such as the Laplace transform and the tail distribution of some random variables, in which many multivariate special functions with symmetric matrix arguments appear, such as gamma, modified Bessel and hypergeometric functions Muirhead [20]. However, the latter case is quite complicated to deal with and, to our knowledge, there are no more precise results on the law of these variables. Nevertheless, in the complex case, hypergeometric functions with Hermitian matrix argument can be expressed as a determinant of a matrix whose entries are one-dimensional hypergeometric functions. In fact, Gross and Richards [10] established the following result:

$$
{}_pF_q(a_1,\ldots,a_p,b_1,\ldots,b_q;X) \\
= \frac{\det(x_i^{n-j}{}_p\mathscr{F}_q(a_1-j+1,\ldots,a_p-j+1,\ldots,b_q-j+1;x_i))}{V(X)}
$$

where $X$ is an $m \times m$ Hermitian matrix, $(x_i)$ are its eigenvalues, ${}_p\mathscr{F}_q$ denotes the standard hypergeometric function with scalar argument, $V(X) = \prod_{i<j}(x_i - x_j)$ is the Vandermonde determinant, and ${}_pF_q$ is the hypergeometric function with Hermitian matrix argument defined by

$$
{}_pF_q(a_1,\ldots,a_p,b_1,\ldots,b_q;X) = \sum_{k\geq 0} \sum_\tau \frac{(a_1)_\tau \cdots (a_p)_\tau}{(b_1)_\tau \cdots (b_q)_\tau} \frac{C_\tau(X)}{k!}
$$

where $\tau = (\tau_1,\ldots,\tau_m)$ is a partition of length at most $m$ and of weight $k$ (i.e. $\tau_1 \geq \tau_2 \geq \cdots \geq \tau_m$, $\sum_i \tau_i = k$), $(a)_\tau$ is the generalized Pochammer symbol and $C_\tau$ is the so-called zonal polynomial. We refer to Macdonald [19] for further details and Lassalle [16, 17] for analogous expressions for multivariate orthogonal polynomials. The determinantal representation above is due to the fact that the zonal polynomial is identified with the (normalized) *Schur functions* defined by:

$$
s_\tau(x_1,\ldots,x_m) = \frac{\det(x_i^{\tau_j+m-j})}{\det(x_i^{m-j})}.
$$

Consequently, one can use integral representations as well as other properties of standard hypergeometric functions to obtain, at least when $m=2$, some results hitherto unknown in the Wishart case. The rest of this paper is organized as follows. In Section 2 we introduce the Laguerre process of integer dimension. In Section 3 we study the behaviour of the eigenvalue process. Then, in Section 4, we define the Laguerre process of positive real dimension. Section 5 is devoted to the absolute-continuity relations, from which we deduce the Laplace transform of the so-called *generalized Hartman–Watson* law as well as the tail distribution of $T_0$, the first hitting time of 0. In Section 6 we focus on the case $m=2$ for which we invert this Laplace transform. In Section 7, we compute the density of $S_0 := 1/(2T_0)$. Section 8 presents our conclusions. Some special functions of interest are discussed in the Appendix.

## 2. Laguerre process of integer index

Let $B$ be a $n \times m$ complex Brownian matrix starting from $B_0$; that is, $B = (B_{ij})$, where the entries $B_{ij}$ are independent complex Brownian motions, so we can write $B = B^1 + \mathrm{i}B^2$, in which $B^1$, $B^2$ are two independent real Brownian matrices. We are interested in the matrix-valued process $X_t := B_t^\star B_t$ which satisfies the following SDE:

$$\mathrm{d}X_t = \mathrm{d}B_t^\star B_t + B_t^\star \, \mathrm{d}B_t + 2nI_m \, \mathrm{d}t. \tag{1}$$

**Definition 1.** $(X_t)_{t \geq 0}$ is called the Laguerre process of size $m$, of dimension $n$ and starting from $X_0 = B_0^\star B_0$, and will be denoted by $L(n, m, X_0)$.

**Remark 1.** For $m = 1$, $(X_t)_{t \geq 0}$ is a squared Bessel process of dimension $2n$, denoted by $\mathrm{BESQ}(2n, X_0)$.

**Remark 2.** Set $X_t = (X_{ij}(t))_{i,j}$. One can easily check that

$$\mathrm{d}X_{ii}(t) = 2\sqrt{X_{ii}(t)} \, \mathrm{d}\gamma_i(t) + 2n \, \mathrm{d}t, \qquad 1 \leq i \leq m,$$

where $(\gamma_i)_{1 \leq i \leq m}$ are independent Brownian motions, thus $X_{ii}$ is a $\mathrm{BESQ}(2n, (X_0)_{ii})$.

**Remark 3.** The equation above implies that

$$\mathrm{d}(\mathrm{tr}(X_t)) = 2\sqrt{\mathrm{tr}(X_t)} \, \mathrm{d}\beta_t + 2nm \, \mathrm{d}t, \tag{2}$$

where $\beta$ is a Brownian motion. Consequently, $(\mathrm{tr}(X_t))_{t \geq 0}$ is a $\mathrm{BESQ}(2nm, \mathrm{tr}(X_0))$ of dimension $2nm$ starting from $\mathrm{tr}(X_0)$. One can also deduce from (1) that, for every $i, j, k, l \in \{1, \ldots, m\}$,

$$\langle \mathrm{d}X_{ij}, \mathrm{d}X_{kl} \rangle_t = 2(X_{il}\delta_{kj} + X_{kj}\delta_{il}) \, \mathrm{d}t,$$

which differs from equation (I-1-5) derived by Bru [2] since, for a complex Brownian motion $\gamma$, we have $\mathrm{d}\langle \gamma, \gamma \rangle_t = 0$ and $\langle \gamma, \overline{\gamma} \rangle_t = 2t$.

Let $H_m, \tilde{H}_m^+$ be respectively the space of $m \times m$ Hermitian matrices and the space of $m \times m$ positive definite Hermitian matrices. On the space of Hermitian matrix-argument functions, we define the matrix-valued differential operators

$$\frac{\partial}{\partial x} := \left(\frac{\partial}{\partial x_{jk}}\right)_{j,k}, \qquad \frac{\partial}{\partial y} := \left(\frac{\partial}{\partial y_{jk}}\right)_{j,k}, \qquad \frac{\partial}{\partial z} := \left(\frac{\partial}{\partial x_{jk}} - \mathrm{i}\frac{\partial}{\partial y_{jk}}\right)_{j,k}.$$

We also define

$$\left(\frac{\partial}{\partial z}\right)^2_{i,j} := \sum_k \frac{\partial^2}{\partial z_{ik}\partial z_{kj}}, \qquad \left(\frac{\partial}{\partial x}\frac{\partial}{\partial y}\right)_{i,j} := \sum_k \frac{\partial^2}{\partial x_{ik}\partial y_{kj}}.$$

**Proposition 1 (Infinitesimal generator).** *Suppose we have functions $f$ satisfying*

$$\frac{\partial f}{\partial x_{ij}} = \frac{\partial f}{\partial x_{ji}}, \qquad \frac{\partial f}{\partial y_{ij}} = -\frac{\partial f}{\partial y_{ji}} \qquad \text{for all } i,j.$$

*Then the infinitesimal generator of a Laguerre process $L(n,m,x)$ is given by*

$$\mathscr{L} = 2n\,\text{tr}\left(\Re\left(\frac{\partial}{\partial z}\right)\right) + 2\left[\text{tr}\left(x\Re\left(\frac{\partial}{\partial z}\right)^2\right) + \text{tr}\left(y\Im\left(\frac{\partial}{\partial z}\right)^2\right)\right], \qquad (3)$$

*where $\partial/\partial z$ is the operator defined above.*

*Remark 4.* Using the fact that $x^{\text{T}} = x$, $y^{\text{T}} = -y$ and $\text{tr}(AB) = \text{tr}(BA) = \text{tr}(B^{\text{T}}A^{\text{T}})$ for any two matrices $A$ and $B$, we can see that

$$\text{tr}\left(y\frac{\partial}{\partial y}\frac{\partial}{\partial x}\right) = \text{tr}\left(\frac{\partial}{\partial x}\frac{\partial}{\partial y}y\right) = \text{tr}\left(y\frac{\partial}{\partial x}\frac{\partial}{\partial y}\right) \quad \Longrightarrow \quad \left(y\Im\left(\frac{\partial}{\partial z}\right)^2\right) = 2\,\text{tr}\left(y\frac{\partial}{\partial x}\frac{\partial}{\partial y}\right).$$

## 3. Eigenvalues of Laguerre process

In this section, we will suppose that $n \geq m$. The following result is due to König and O'Connell [15], Katori and Tanemura [14] and Bru [1] in the real case:

**Theorem 1.** *Let $\lambda_1(t), \ldots, \lambda_m(t)$ denote the eigenvalues of $X_t$. Suppose that at time $t = 0$ all the eigenvalues are distinct. Then the eigenvalue process $(\lambda_1(t), \ldots, \lambda_m(t))$ satisfies the following stochastic differential system:*

$$\mathrm{d}\lambda_i(t) = 2\sqrt{\lambda_i(t)}\,\mathrm{d}\beta_i(t) + 2\left[n + \sum_{k\neq i}\frac{\lambda_i(t) + \lambda_k(t)}{\lambda_i(t) - \lambda_k(t)}\right]\mathrm{d}t, \qquad 1 \leq i \leq m,\ t < \tau,$$

*where the $(\beta_i)_{1\leq i\leq m}$ are independent Brownian motions and $\tau$ is defined by $\tau := \inf\{t, \lambda_i(t) = \lambda_j(t)$ for some $(i,j)\}$.*

*Remark 5.* With the help of the SDE satisfied by the eigenvalues, we can compute the ones satisfied by both processes $(\text{tr}(X))$ and $(\det(X))$: the first has already been done. For the second, we find that for $t < T_0 := \inf\{t, \det(X_t) = 0\}$ and for $r \in \mathbb{R}$,

$$\mathrm{d}(\det(X_t)) = 2\det(X_t)\sqrt{\text{tr}(X_t^{-1})}\,\mathrm{d}\nu_t + 2(n - m + 1)\det(X_t)\,\text{tr}(X_t^{-1})\,\mathrm{d}t,$$

$$\mathrm{d}(\log(\det(X_t))) = 2\sqrt{\text{tr}(X_t^{-1})}\,\mathrm{d}\nu_t + 2(n - m)\,\text{tr}(X_t^{-1})\,\mathrm{d}t,$$

$$\mathrm{d}(\det(X_t)^r) = 2r(\det(X_t))^r\sqrt{\text{tr}(X_t^{-1})}\,\mathrm{d}\nu_t + 2r(n - m + r)(\det(X_t))^r\,\text{tr}(X_t^{-1})\,\mathrm{d}t,$$

so we can see that for $n = m$, $\log(\det(X))$ is a local martingale and so is $(\det(X))^{m-n}$.

**Lemma 1.** *Take $X_0 \in \tilde{H}_m^+$. Then for $n \geq m$, $X_t \in \tilde{H}_m^+$.*

**Proof.** In fact, this result is a direct consequence of the fact that, for $n = m$, $\log \det(X)$ is a local martingale, and so is $(\det(X))^{m-n}$. Hence, for $n \geq m$, these two continuous processes tend to infinity when $t \to T_0$, which is possible only if $T_0 = \infty$, because every continuous local martingale is a time-changed Brownian motion. □

**Corollary 1.** *If $\lambda_1(0) > \cdots > \lambda_m(0)$, then the process $U$ defined by*

$$U(t) = \frac{1}{\prod_{i<j}(\lambda_i(t) - \lambda_j(t))}, \qquad t < \tau,$$

*is a local martingale.*

**Proof.** We could follow the proof given by Bru [1] or make straightforward computations using the derivatives of the Vandermonde function. But we prefer to use a result from König, W. and O'Connell [15]: for $n \geq m$, the eigenvalue process is the $V$-transform (in the Doob sense) of the process obtained from $m$ independent BESQ$(2(n-m+1))$. Thus, if $G$ and $\hat{G}$ denote respectively the infinitesimal generators of these two processes, then, $G(h) = 0$ and, for all $C^2$ function $f$,

$$\hat{G}(f) = \frac{1}{V}G(Vf) \quad \Longrightarrow \quad \hat{G}(U) = \frac{1}{V}G(\mathbf{1}) = 0. \qquad \square$$

**Corollary 2.** *If, at time $t = 0$, the eigenvalues of $X$ are distinct, then they will never collide; that is, $\tau = \infty$ almost surely.*

**Proof.** This result follows from the fact that the continuous process $U$ tends to infinity when $t \to \tau$, which is possible only if $\tau = \infty$ almost surely (We use the same argument as before.) □

The proof of the following result is similar to the one derived by Bru [2] in the real case:

**Proposition 2 (Additivity property).** *If $(X_t)_{t \geq 0}$ and $(Y_t)_{t \geq 0}$ are two independent Laguerre processes $L(n, m, X_0)$ and $L(p, m, Y_0)$ respectively, then the process $(X_t + Y_t)_{t \geq 0}$ is a Laguerre process $L(n + p, m, X_0 + Y_0)$.*

We now introduce the Laguerre processes of non-integer dimensions $\delta$.

## 4. Laguerre processes with non-integer dimensions

Let $X$ be a Laguerre process $L(n, m, X_0)$ with $n \geq m$. If $X_0 \in \tilde{H}_m^+$, and if $\sqrt{X_t}$ stands for the symmetric matrix square root of $X_t$, it is easy to show that the matrix $O$ defined by $O_t := \sqrt{X_t}^{-1} B_t^\star$, where $X_t = B_t^\star B_t$, satisfies $O^\star O = OO^\star = I_m$. Thus,

$$\mathrm{d}\gamma_t = O_t\, \mathrm{d}B_t = \sqrt{X_t}^{-1} B_t^\star\, \mathrm{d}B_t$$

is an $m \times m$ complex Brownian matrix. Replacing this expression in (1), one obtains

$$\mathrm{d}X_t = \sqrt{X_t}\, \mathrm{d}\gamma_t + \mathrm{d}\gamma_t^\star \sqrt{X_t} + 2nI_m\, \mathrm{d}t.$$

**Theorem 2.** *If $(B_t)$ is an $m \times m$ complex Brownian matrix, then for every $X_0 \in \tilde{H}_m^+$ and for all $\delta \geq m$, the SDE*

$$\mathrm{d}X_t = \sqrt{X_t}\, \mathrm{d}B_t + \mathrm{d}B_t^\star \sqrt{X_t} + 2\delta I_m\, \mathrm{d}t \tag{4}$$

*has a unique strong solution in $\tilde{H}_m^+$. Furthermore, if the eigenvalues are distinct at time $t = 0$, then they satisfy the stochastic differential system*

$$\mathrm{d}\lambda_i(t) = 2\sqrt{\lambda_i(t)}\, \mathrm{d}\beta_i(t) + 2\left[\delta + \sum_{k \neq i} \frac{\lambda_i(t) + \lambda_k(t)}{\lambda_i(t) - \lambda_k(t)}\right] \mathrm{d}t, \qquad 1 \leq i \leq m,$$

*where the $(\beta_i)_{1 \leq i \leq m}$ are independent Brownian motions.*

**Proof.** The proof of the second part of the theorem is the same as before, with $\delta$ instead of $n$. All that remains is to prove the first part. Note first that $(\det(X_t)), (\log \det(X_t))$ and $(\det(X_t)^r)$ satisfy the same SDE with $\delta$ instead of $n$. Hence, arguing as before, we can see that $T_0 = \infty$ almost surely. On the other hand, the map $a \mapsto a^{1/2}$ is analytic in $\tilde{H}_m^+$ (see Rogers and Williams [22], page 134), so the SDE has a unique strong solution for all $t \geq 0$. □

**Definition 2.** *Such a process is called the* Laguerre process *of dimension $\delta$, size $m$ and initial state $X_0$. It will be denoted by $L(\delta, m, X_0)$.*

**Remark 6.** Any process $(X_t)_{t \geq 0}$ solution of (4) is a diffusion whose infinitesimal generator is given by

$$\mathscr{L} = 2\delta\, \mathrm{tr}\left(\Re\left(\frac{\partial}{\partial z}\right)\right) + 2\left[\mathrm{tr}\left(x\Re\left(\frac{\partial}{\partial z}\right)^2\right) + \mathrm{tr}\left(y\Im\left(\frac{\partial}{\partial z}\right)^2\right)\right].$$

**Remark 7.** A simple computation shows that

$$\mathrm{d}\langle X_{ij}, X_{kl}\rangle_t = 2(X_{il}(t)\delta_{kj} + X_{kj}(t)\delta_{il})\, \mathrm{d}t, \qquad \text{for all } i, j, k, l \in \{1, \ldots, m\}.$$

We now focus on both existence and uniqueness when $\delta > m - 1$ and $X_0 \in H_m^+$ – see Bru [2] for the real case.

If $X$ is a Hermitian matrix, let $X^+$ be the Hermitian matrix $\max(X, 0)$. If we denote by $(\lambda_i)$ the eigenvalues of $X$, then $(\lambda_i^+ = \max(\lambda_i, 0))$ are those of $X^+$.

**Theorem 3.** *For all $\delta \in \mathbb{R}_+$ and $X_0 = x \in H_m$, the SDE*

$$\mathrm{d}X_t = \sqrt{X_t^+}\,\mathrm{d}B_t + \mathrm{d}B_t^\star\sqrt{X_t^+} + 2\delta I_m\,\mathrm{d}t \tag{5}$$

*has a solution in $H_m$.*

**Proof.** The mapping $a \mapsto \sqrt{a^+}$ is continuous on $H_m$. Hence, $X$ exists up to its explosion time (Ikeda and Watanabe [11], page Theorem 2.3). Furthermore, from

$$\|\sqrt{X^+}\|^2 + \|\delta I\|^2 \leq \delta^2 + \|X\|^2 \leq C(1 + \|X\|^2),$$

we can deduce that this explosion time is infinite almost surely (Ikeda and Watanabe [11], page Theorem 2.4). $\square$

**Proposition 3.** *If $\lambda_1(0) > \cdots > \lambda_m(0) \geq 0$, then, for all $t < S := \inf\{t, \lambda_i = \lambda_j \text{ for some } (i,j)\}$, the eigenvalues of $X^+$ satisfy the following differential system:*

$$\mathrm{d}\lambda_i(t) = 2\sqrt{\lambda_i^+(t)}\,\mathrm{d}\nu_i(t) + 2\left(\delta + \sum_{k \neq i} \frac{\lambda_i^+(t) + \lambda_k^+(t)}{\lambda_i(t) - \lambda_k(t)}\right)\mathrm{d}t, \qquad 1 \leq i \leq m.$$

**Proof.** This differential system can be shown in the same way as in Theorem 1, using

$$\langle \mathrm{d}X_{ij}, \mathrm{d}X_{kl}\rangle_t = 2(X_{il}^+(t)\delta_{kj} + X_{kj}^+(t)\delta_{il})\,\mathrm{d}t, \qquad \text{for all } i, j, k, l \in \{1, \ldots, m\}. \quad \square$$

**Proposition 4.** *If $\lambda_1(0) > \cdots > \lambda_m(0) \geq 0$, then, for all $\delta > m - 1$, $t > 0$, $\lambda_m(t) \geq 0$.*

**Proof.** First, we note that $S = \infty$ almost surely. Indeed, one can easily show that the process $U$ defined by

$$U(\lambda_1(t), \ldots, \lambda_m(t)) = \frac{1}{\prod_{i<j}(\lambda_i(t) - \lambda_j(t))}$$

is a local martingale. For the proof, we proceed along the lines of Bru [1]. $\square$

**Theorem 4.** *If $\lambda_1(0) > \cdots > \lambda_m(0) \geq 0$, then, for all $\delta > m-1$, (4) has a unique solution in $H_m^+$ in the sense of probability law.*

**Proof.** By Proposition 4, the solution of the SDE (5) remains positive for all $t > 0$, thus it is a solution of (4). $\square$

**Theorem 5.** *Let $H_m^+$ be the space of positive Hermitian matrices. Then, whenever the SDE (4) has a solution in $H_m^+$, for fixed $t$, its distribution is given by its Laplace transform*

$$\mathbb{E}_{X_0}(\exp(-\operatorname{tr} uX_t)) = (\det(I_m + 2tu))^{-\delta} \exp(-\operatorname{tr}(X_0(I_m + 2tu)^{-1}u)), \tag{6}$$

*for all $u$ in $H_m^+$.*

**Proof.** For $s \in H_m^+$, let $g(t,s) = \Delta_t^{-\delta} \exp(-V(t,s))$, where

$$\Delta_t = \det(I_m + 2ut), \qquad W_t = (I_m + 2ut)^{-1} u, \qquad V(t,s) = \operatorname{tr}(sW_t).$$

First, note that $W \in H_m$. To proceed, we need the following lemma:

**Lemma 2.** *The function $g$ satisfies the heat equation: $\partial g / \partial t = \mathscr{L} g$, where $\mathscr{L}$ is the infinitesimal generator of $X$.*

**Proof.** If we write $s = x + iy$, then, using the fact that $x$ is symmetric, $y$ is skew-symmetric and $W$ is Hermitian, we can see that $\operatorname{tr}(sW_t) = \operatorname{tr}(xM + iyN)$, where

$$M = \frac{W + \overline{W}}{2}, \qquad N = \frac{W - \overline{W}}{2}.$$

Observing that $M^\mathrm{T} = M$ and $N^\mathrm{T} = -N$, we can deduce that $g$ satisfies the conditions of Proposition 1. Furthermore,

$$\frac{\partial g}{\partial t} = -g(2\delta \operatorname{tr}(W_t) - 2 \operatorname{tr}(sW_t^2))$$

$$\operatorname{tr}\left( y \left( \frac{\partial^2 g}{\partial x \partial y} + \frac{\partial^2 g}{\partial y \partial x} \right) \right) = -ig \operatorname{tr}(yW^2),$$

$$\operatorname{tr}\left( x \left( \frac{\partial^2 g}{\partial x^2} - \frac{\partial^2 g}{\partial y^2} \right) \right) = g \operatorname{tr}(x(M^2 + N^2)) = g \operatorname{tr}(xW^2).$$

Finally, noting that $\operatorname{tr}(M) = \operatorname{tr}(W)$, we obtain the equality.

Now, we consider the process $(Z(t, X_t))$ defined by $Z(t, X_t) = g(t_1 - t, X_t)$ for all $t \leq t_1$ for fixed $t_1$. From the lemma, we deduce that $Z$ is a bounded local martingale and thus is a martingale. So, the result follows from a simple application of the optional stopping theorem. □

**Corollary 3.** *Let $(X_t)_{t \geq 0}$ be a Laguerre process $L(\delta, m, x)$ where $x \in \tilde{H}_m^+$. For $\delta > m - 1$, its semi-group is given by the following density:*

$$p_t^\delta(x, y) = \frac{1}{(2t)^{m\delta} \Gamma_m(\delta)} \exp\left( -\frac{1}{2t} \operatorname{tr}(x + y) \right) (\det y)^{\delta - m} {}_0F_1\left( \delta; \frac{xy}{4t^2} \right) \mathbf{1}_{\{y > 0\}}$$

*with respect to Lebesgue measure $\mathrm{d}y = \prod_{p \leq q} \mathrm{d}y_{pq}^1 \prod_{p < q} \mathrm{d}y_{pq}^2$ where $y = y^1 + iy^2$ and ${}_0F_1$ is a hypergeometric function with matrix argument (Chikuze [4]; Gross and Richard [10]).*

**Proof of Theorem 5.** In fact, this result can be easily deduced from the case where $\delta = n$ is integer, since, in this case, $X_t$ is a non-central complex Wishart variable $W(n, 2tI_m, x)$ James [12] with density given by

$$f_t(x, y) = \frac{1}{(2t)^{mn} \Gamma_m(n)} \exp\left( -\frac{1}{2t} \operatorname{tr}(x + y) \right) (\det y)^{n-m} {}_0F_1\left( n; \frac{xy}{4t^2} \right) \mathbf{1}_{\{y > 0\}}$$

with respect to $dy$. Hence, taking $\delta$ instead of $n$ and denoting by $W_t$ this new variable (starting from $x$), we can see that, using $|y|$ to denote $\det(y)$,

$$\begin{aligned}
E_x(e^{-\operatorname{tr} uW_t}) &= \frac{1}{(2t)^{m\delta}\Gamma_m(\delta)} \exp\left(-\frac{\operatorname{tr} x}{2t}\right) \\
&\quad \times \int_{y>0} \exp\left(-\frac{1}{2t}\operatorname{tr}((I+2ut)y)\right)|y|^{\delta-m} {}_0F_1\left(\delta; \frac{xy}{4t^2}\right) dy \\
&= \frac{2t^{m\delta}|x|^{-\delta}}{\Gamma_m(\delta)} \exp\left(-\frac{\operatorname{tr} x}{2t}\right) \\
&\quad \times \int_{z>0} \exp(-2t\operatorname{tr}(x^{-1/2}(I+2ut)x^{-1/2}z))|z|^{\delta-m} {}_0F_1(\delta; z) dz \\
&= \exp\left(-\frac{\operatorname{tr} x}{2t}\right)|I+2ut|^{-\delta}\exp\left(\operatorname{tr}\left(\frac{x}{2t}(I+2ut)^{-1}\right)\right) \\
&= |I+2ut|^{-\delta}\exp\left(-\frac{1}{2t}\operatorname{tr}(x(I+2ut)^{-1}(I+2ut-I))\right) \\
&= |I+2ut|^{-\delta}\exp\left(-\operatorname{tr}(x(I+2ut)^{-1}u)\right),
\end{aligned}$$

which is equal to (6). □

**Remark 8.** In the last proof, we used the change of variables $z = x^{1/2}yx^{1/2}$ which gives $dz = |x|^m dy$. For the second integral, see Faraut and Korányi [8], Proposition XV.1.3.

**Remark 9.** The expression for the semi-group extends continuously to the degenerate case, namely:

$$p_t^\delta(0_m, y) = \frac{1}{(2t)^{m\delta}\Gamma_m(\delta)} \exp\left(-\frac{\operatorname{tr}(y)}{2t}\right)(\det y)^{\delta-m}\mathbf{1}_{\{y>0\}},$$

where $0_m$ denotes the null matrix.

**Corollary 4.** *For $\delta > m-1$, the semi-group of eigenvalue process is given by*

$$q_t(x,y) = \frac{V(y)}{V(x)}\det\left(\frac{1}{2t}\left(\frac{y_j}{x_i}\right)^{\nu/2} e^{-(x_i+y_j)/(2t)} I_\nu\left(\frac{\sqrt{x_iy_j}}{t}\right)\right),$$

*where $x = (x_1, \ldots, x_m), y = (y_1, \ldots, y_m)$ so that $x_1 > \cdots > x_m > 0, y_1 > \cdots > y_m > 0$, $\delta = m + \nu$ such that $\nu > -1$ and $I_\nu$ denotes the modified Bessel function Lebedev [18].*

**Proof.** The expression for the semi-group can be computed using Karlin and MacGregor [13] formula since, for $\delta > m-1$, the eigenvalue process is the $h$-transform of the process consisting of $m$ independent BESQ$(2(\delta - m + 1))$ conditioned never to collide, as stated by König and O'Connell [15]. Another proof is given by Péché ([21], page 68). Here, we will deduce the expression for $q_t(x,y)$ from $p_t(x,y)$ following Muirhead [20], namely, by projection on the unitary group: we will use the Weyl integration formula, then give a determinantal representation of hypergeometric functions of two matrix arguments. First, we state the Weyl integration formula Faraut [7] in the complex case: for any Borel function $f$,

$$\int_{H_m} f(A)\,\mathrm{d}A = C_m \int_{U(m)} \int_{\mathbb{R}^m} f(uau^*)\alpha(\mathrm{d}u)(V(a))^2\,\mathrm{d}a_1 \cdots \mathrm{d}a_m,$$

where $C_m = \pi^{m(m-1)}/\Gamma_m(m)$, $U(m)$ is the unitary group, $\alpha$ is the normalized Haar measure on $U(m)$, $a = \mathrm{diag}(a_i)$ and $A = uau^*$. Hence, the semi-group of the eigenvalue process is given by James [12]:

$$\begin{aligned}
q_t(x,y) &= C_m(V(y)^2) \int_{U(m)} p_t(\tilde{x}, u\tilde{y}u^*)\alpha(\mathrm{d}u) \\
&= \frac{C_m(V(y)^2)}{(2t)^{m\delta}\Gamma_m(\delta)} \prod_{i,j=1}^m \mathrm{e}^{-(x_i+y_j)/2t} \left(\prod_{i=1}^m y_j\right)^{\delta-m} \int_{U(m)} {}_0F_1\left(\delta; \frac{\tilde{x}u\tilde{y}u^*}{4t^2}\right)\alpha(\mathrm{d}u) \\
&= \frac{\pi^{m(m-1)}(V(y)^2)}{(2t)^{m(m+\nu)}\Gamma_m(m)\Gamma_m(m+\nu)} \prod_{i,j=1}^m \mathrm{e}^{-(x_i+y_j)/2t} \left(\prod_{i=1}^m y_j\right)^{\nu} {}_0F_1\left(m+\nu; \frac{\tilde{x}}{4t^2}; \tilde{y}\right),
\end{aligned}$$

where $\tilde{y} = \mathrm{diag}(y_j)$, $x$ is a positive definite matrix whose eigenvalues are $x_1, \ldots, x_m$, ${}_0F_1$ in the third line is a hypergeometric function with two matrix arguments Gross and Richards [10] and $\delta = m + \nu$, $\nu > -1$. Next, we need a lemma.

**Lemma 3.** *Let $B, C \in H_m$ and let $(b_i), (c_i)$ be their respective eigenvalues. Then*

$$\begin{aligned}
{}_pF_q((m+\mu_i)_{1\leq i\leq p}, (m+\phi_j)_{1\leq j\leq q}; B, C) \\
= \pi^{m(m-1)/2(p-q-1)}\Gamma_m(m) \\
\times \prod_{i=1}^p \frac{(\Gamma(\mu_i+1))^m}{\Gamma_m(m+\mu_i)} \\
\times \prod_{j=1}^q \frac{\Gamma_m(m+\phi_j)}{(\Gamma(\phi_j+1))^m} \frac{\det({}_p\mathscr{F}_q((\mu_i+1)_{1\leq i\leq p}, (1+\phi_j)_{1\leq j\leq q}; b_lc_f)_{l,f})}{V(B)V(C)}
\end{aligned}$$

*for all $\mu_i, \phi_j > -1, 1 \leq i \leq p, 1 \leq j \leq q$.*

**Proof.** Recall that the hypergeometric function of two matrix arguments is given by the series

$$_pF_q((a_i)_{1\leq i\leq p}, (e_j)_{1\leq j\leq q}; B, C) = \sum_{k=0}^{\infty}\sum_{\tau} \frac{\prod_{i=1}^{p}(a_i)_\tau}{\prod_{j=1}^{q}(e_j)_\tau} \frac{C_\tau(B)C_\tau(C)}{C_\tau(I)k!}.$$

It is well known that

$$C_\tau(B) = \frac{k! d_\tau}{(m)_\tau} s_\tau(b_1, \ldots, b_m),$$

where $s_\tau$ is the Schur function and $d_\tau = s_\tau(1, \ldots, 1)$ is the representation trace or degree (Gross and Richards [10] and Faraut [7]). Substituting in the series gives

$$_pF_q((m+\mu_i)_{1\leq i\leq p}, (m+\phi_j)_{1\leq j\leq q}; B, C) = \sum_{k=0}^{\infty}\sum_{\tau} \frac{\prod_{i=1}^{p}(m+\mu_i)_\tau}{\prod_{j=1}^{q}(m+\phi_j)_\tau} \frac{s_\tau(B)s_\tau(C)}{(m)_\tau}.$$

We now write

$$(m+\mu_i)_\tau = \prod_{r=1}^{m} \frac{\Gamma(\mu_i+m+k_r-r+1)}{\Gamma(\mu_i+m-r+1)} = \prod_{r=1}^{m} \frac{\Gamma(\mu_i+1+k_r+\delta_r)}{\Gamma(\mu_i+m-r+1)}$$

$$= \pi^{m(m-1)/2} \frac{(\Gamma(\mu_i+1))^m}{\Gamma_m(m+\mu_i)} \prod_{r=1}^{m}(\mu_i+1)_{k_r+\delta_r}, \qquad \delta_r = m-r.$$

Doing the same for each $(m+\phi_j)_\tau$ and for $(m)_\tau$, we can see that

$$_pF_q((m+\mu_i)_{1\leq i\leq p}, (m+\phi_j)_{1\leq j\leq q}; B, C)$$
$$= \pi^\beta \Gamma_m(m) \prod_{i=1}^{p} \frac{(\Gamma(\mu_i+1))^m}{\Gamma_m(m+\mu_i)} \prod_{j=1}^{q} \frac{(\Gamma(\phi_j+1))^m}{\Gamma_m(m+\phi_j)}$$
$$\times \sum_{k=0}^{\infty}\sum_{\tau}\prod_{r=1}^{m}\left(\frac{\prod_{i=1}^{p}(\mu_i+1)_{k_r+\delta_r}}{\prod_{j=1}^{q}(\phi_j+1)_{k_r+\delta_r}}\right) \frac{s_\tau(B)s_\tau(C)}{\prod_{r=1}^{m}(1)_{k_r+\delta_r}}$$

where

$$\beta = \left(\frac{m(m-1)}{2}\right)(p-q-1).$$

To obtain the desired result, we use the Hua formula Faraut [7]:

**Lemma 4.** *Given an entire function* $f(z) = \sum_{k=0}^{\infty} e_k z^k$, *we have*

$$\frac{\det(f(b_i c_j))_{i,j}}{V(B)V(C)} = \sum_{k=0}^{\infty}\sum_{\tau}\left(\prod_{r=1}^{m} e_{k_r+\delta_r}\right) \frac{s_\tau(B)s_\tau(C)}{s_\tau(I_m)}.$$

Thus, we obtain

$$_pF_q((m+\mu_i)_{1\leq i\leq p},(m+\phi_j)_{1\leq j\leq q};B,C)$$
$$=\pi^{m(m-1)(p-q-1)/2}\Gamma_m(m)\prod_{i=1}^{p}\frac{\Gamma(\mu_i+1)}{\Gamma_m(m+\mu_i)}$$
$$\times\prod_{j=1}^{q}\frac{\Gamma_m(m+\phi_j)}{\Gamma(\phi_j+1)}\frac{\det(\sum_{k=0}^{\infty}\frac{\prod_{i=1}^{p}(\mu_i+1)_k}{\prod_{j=1}^{q}(\phi_j+1)_k}\frac{(b_lc_p)^k}{k!})_{l,p}}{V(B)V(C)}$$
$$=\pi^{m(m-1)(p-q-1)/2}\Gamma_m(m)\prod_{i=1}^{p}\frac{\Gamma(\mu_i+1)}{\Gamma_m(m+\mu_i)}\prod_{j=1}^{q}\frac{\Gamma_m(m+\phi_j)}{\Gamma(\phi_j+1)}$$
$$\times\frac{\det(_p\mathscr{F}_q((\mu_i+1)_{1\leq i\leq p},(1+\phi_j)_{1\leq j\leq q};b_lc_f)_{l,f})}{V(B)V(C)}.\qquad\square$$

For $p=0$ and $q\geq 1$,

$$_0F_q((m+\phi_j)_{1\leq j\leq q};B,C)$$
$$=\pi^{-m(m-1)(q+1)/2}\Gamma_m(m)\prod_{j=1}^{q}\frac{\Gamma_m(m+\phi_j)}{(\Gamma(\phi_j+1))^m}\frac{\det(_0\mathscr{F}_q((1+\phi_j)_{1\leq j\leq q};b_lc_f)_{l,f}}{V(B)V(C)},$$

and similarly,

$$_0F_0(B,C)=\frac{\Gamma_m(m)}{\pi^{m(m-1)/2}}\frac{\det(e^{b_lc_f})_{l,f}}{V(B)V(C)}$$

which can be viewed as a Harish–Chandra formula for the Itzykson–Zuber integral Collins [5]. We now proceed to the end of the proof. Taking $p=0$, $q=1$, $B=\tilde{x}/4t^2$, $C=\tilde{y}$, we obtain

$$_0F_1\left(m+\nu;\frac{\tilde{x}}{4t^2};\tilde{y}\right)=\frac{(4t^2)^{m(m-1)/2}\Gamma_m(m+\nu)\Gamma_m(m)}{\pi^{m(m-1)}(\Gamma(\nu+1))^m}\frac{\det(_0\mathscr{F}_1((\nu+1);x_iy_j/4t^2))}{V(x)V(y)}.$$

The expression for $q_t(x,y)$ follows from a simple computation and from the fact that

$$\frac{_0\mathscr{F}_1((\nu+1);x_iy_j/4t^2)_{i,j}}{\Gamma(\nu+1)}=\left(\frac{2t}{\sqrt{x_iy_j}}\right)^{\nu}I_{\nu}\left(\frac{\sqrt{x_iy_j}}{t}\right).\qquad\square$$

**Proposition 5.** *The measure defined by $\rho(\mathrm{d}x)=(\det(x))^{\delta-m}\,\mathrm{d}x$ on $\tilde{H}_m^+$ is invariant under the semi-group, that is, $\rho P_t=\rho$.*

**Proof.** Denote by $P_t$ the semi-group of the Laguerre process $L(\delta, m, x)$ for $\delta > m - 1$. Then we have to show that

$$\int_{x>0} P_t f(x) \rho(\mathrm{d}x) = \int_{y>0} f(y) \rho(\mathrm{d}y), \qquad f \in C_0(\tilde{H}_m^+).$$

This follows by a similar computation and the same arguments as in the proof of Corollary 3. □

**Remark 10.** For Wishart processes, it is easy to see that

$$\mu(\mathrm{d}x) := (\det(x))^{\delta/2 - (m+1)/2} \mathbf{1}_{\{x>0\}} \, \mathrm{d}x$$

is invariant under the semi-group.

## 5. Girsanov formula and absolute-continuity relations

The index $\nu > -1$ of a $L(\delta, m, x)$ is defined by $\nu = \delta - m$. In this section, we will proceed along the same lines as Donati-Martin *et al.* [6] to derive absolute-continuity relations between different indices.

### 5.1. Positive indices

Take a matrix-valued Hermitian predictable process $H$. Let $Q_x^\delta$ be the probability law of $L(\delta, m, x)$ for $\delta > m - 1$ and $x \in \tilde{H}_m^+$. Define

$$L_t = \int_0^t \frac{\operatorname{tr}(H_s \, \mathrm{d}B_s + \overline{H_s} \, \overline{\mathrm{d}B_s})}{2},$$

$$\Phi_t = \exp\left(L_t - \frac{1}{2} \int_0^t \operatorname{tr}(H_s^2) \, \mathrm{d}s\right),$$

where $B$ is a complex Brownian matrix under $Q_x^\delta$. We can easily see that the process $\beta$ defined by $\beta_t = B_t - \int_0^t H_s \, \mathrm{d}s$ is a Brownian matrix under the probability

$$\mathbb{P}_x^H|_{\mathscr{F}_t} := \Phi_t \cdot Q_x^\delta|_{\mathscr{F}_t}.$$

Furthermore, $(X_t)_{t \geq 0}$ is a solution of

$$\mathrm{d}X_t = \sqrt{X_t} \, \mathrm{d}\beta_t + \mathrm{d}\beta_t^\star \sqrt{X_t} + (\sqrt{X_t} H_t + H_t \sqrt{X_t} + 2\delta I_m) \, \mathrm{d}t. \qquad (7)$$

For $H_t = \nu \sqrt{X_t}^{-1}$, (7) becomes

$$\mathrm{d}X_t = \sqrt{X_t} \, \mathrm{d}\beta_t + \mathrm{d}\beta_t^\star \sqrt{X_t} + 2(\delta + \nu) I_m \, \mathrm{d}t,$$

so that $(X_t)_{t \geq 0}$ is $L(\delta + \nu, m, x)$ under $\mathbb{P}_x^H$. Thus, we have proved the following theorem.

**Theorem 6.** *For $\delta > m - 1$,*

$$Q_x^{\delta+\nu}|_{\mathscr{F}_t} = \exp\left(\frac{\nu}{2}\int_0^t \mathrm{tr}\left(\sqrt{X_s}^{-1}\mathrm{d}B_s + \overline{\sqrt{X_s}^{-1}\mathrm{d}B_s}\right) - \frac{\nu^2}{2}\int_0^t \mathrm{tr}(X_s^{-1})\,\mathrm{d}s\right) \cdot Q_x^{\delta}|_{\mathscr{F}_t}.$$
(8)

**Proposition 6.**

$$Q_x^{m+\nu}|_{\mathscr{F}_t} = \left(\frac{\det(X_t)}{\det(x)}\right)^{\nu/2} \exp\left(-\frac{\nu^2}{2}\int_0^t \mathrm{tr}(X_s^{-1})\,\mathrm{d}s\right) \cdot Q_x^{m}|_{\mathscr{F}_t}.$$
(9)

**Proof.** We know that $\nabla_u(\det(u)) = \det(u)u^{-1}$, hence, $\nabla_u(\log(\det(u))) = u^{-1}$. Then, using the fact that, for $\delta = m$, $(\log(\det(X_t)))$ is a local martingale, we obtain from the Itô formula that

$$\log(\det(X_t)) = \log(\det(X_0)) + \int_0^t \mathrm{tr}(X_s^{-1}(\sqrt{X_s}\mathrm{d}B_s + \mathrm{d}B_s^\star\sqrt{X_s}))$$

$$= \log(\det(X_0)) + \int_0^t \mathrm{tr}\left(\sqrt{X_s}^{-1}\mathrm{d}B_s + \overline{\sqrt{X_s}^{-1}\mathrm{d}B_s}\right). \qquad \square$$

From (9), it follows that:

**Corollary 5.**

$$Q_x^m\left(\exp\left(-\frac{\nu^2}{2}\int_0^t \mathrm{tr}(X_s^{-1})\,\mathrm{d}s\right)\bigg|X_t = y\right)$$

$$= \frac{\det(y)}{\det(x)}^{-\nu/2}\frac{p_t^{m+\nu}(x,y)}{p_t^m(x,y)} = \frac{\Gamma_m(m)}{\Gamma_m(m+\nu)}(\det(z))^{\nu/2}\frac{{}_0F_1(m+\nu,z)}{{}_0F_1(m,z)} := \frac{\tilde{I}_\nu(z)}{\tilde{I}_0(z)},$$

*where $z = xy/4t^2$.*

We now state the following asymptotic result:

**Corollary 6.** *Let $X$ be a Laguerre process $L(m, m, x)$. Then, as $t \to \infty$,*

$$\frac{4}{(m\log t)^2}\int_0^t \mathrm{tr}(X_s)^{-1}\,\mathrm{d}s \xrightarrow{\mathscr{L}} T_1(\beta),$$

*where $T_1$ is the first hitting time of 1 by a standard Brownian motion $\beta$.*

**Proof.** From (9), we deduce that

$$Q_x^m\left(\exp\left(-\frac{2\nu^2}{(m\log t)^2}\int_0^t \mathrm{tr}(X_s^{-1})\,\mathrm{d}s\right)\bigg|X_t = ty\right)$$

$$= \frac{\Gamma_m(m)}{\Gamma_m(m+2\nu/m\log t)}(\det(xy/4t))^{\nu/m\log t}\frac{{}_0F_1(m+2\nu/m\log t,xy/4t^2)}{{}_0F_1(m,xy/4t^2)}.$$

Noting that $(t^m)^{-\nu/m\log t}=\mathrm{e}^{-\nu}$, and since both hypergeometric functions converge to 1 as $t\to\infty$, we obtain

$$Q_x^m\left(\exp\left(-\frac{2\nu^2}{(m\log t)^2}\int_0^t \mathrm{tr}(X_s^{-1})\,\mathrm{d}s\Big|X_t=ty\right)\right) \stackrel{t\to\infty}{\longrightarrow} \mathrm{e}^{-\nu}.$$

Then, since

$$\lim_{t\to\infty} t^{m^2} p_t^m(x,2y) = \lim_{t\to\infty} \frac{\mathrm{e}^{-\mathrm{tr}(x)/2t}}{\Gamma_m(m)}\mathrm{e}^{-\mathrm{tr}(y)}\,{}_0F_1\left(m,\frac{xy}{2t}\right) = \frac{\mathrm{e}^{-\mathrm{tr}(y)}}{\Gamma_m(m)},$$

we obtain

$$Q_x^m\left(\exp\left(-\frac{2\nu^2}{(m\log t)^2}\int_0^t \mathrm{tr}(X_s^{-1})\,\mathrm{d}s\right)\right)$$
$$= \int_{y>0} Q_x^m\left(\exp\left(-\frac{2\nu^2}{(m\log t)^2}\int_0^t \mathrm{tr}(X_s^{-1})\,\mathrm{d}s\right)\Big|X_t=y\right) p_t^m(x,y)\,\mathrm{d}y$$
$$= \int_{y>0} Q_x^m\left(\exp\left(-\frac{2\nu^2}{(m\log t)^2}\int_0^t \mathrm{tr}(X_s^{-1})\,\mathrm{d}s\right)\Big|X_t=ty\right) t^{m^2} p_t^m(x,ty)\,\mathrm{d}y \stackrel{t\to\infty}{\longrightarrow} \mathrm{e}^{-\nu}$$

by the dominated convergence theorem. □

## 5.2. Negative indices

Take $0 < a \leq \det(x)$. The same computation as in Section 5.1 with $H_t = -\nu\sqrt{X_t}^{-1}$, $0 < \nu < 1$, shows that

$$Q_x^{m-\nu}|_{\mathscr{F}_{t\wedge T_a}} = \left(\frac{\det(x)}{\det(X_{t\wedge T_a})}\right)^{\nu/2} \exp\left(-\frac{\nu^2}{2}\int_0^{t\wedge T_a} \mathrm{tr}(X_s^{-1})\,\mathrm{d}s\right) Q_x^m|_{\mathscr{F}_{t\wedge T_a}}$$

where $T_a := \inf\{t, \det(X_t) = a\}$. Letting $a \to 0$ and using the fact that $T_0 = \infty$ almost surely under $Q_x^m$, we obtain

$$Q_x^{m-\nu}|_{\mathscr{F}_{t\wedge T_0}} = \left(\frac{\det(x)}{\det(X_t)}\right)^{\nu/2} \exp\left(\frac{\nu^2}{2}\int_0^t \mathrm{tr}(X_s^{-1})\,\mathrm{d}s\right) Q_x^m|_{\mathscr{F}_t}$$
$$= \left(\frac{\det(x)}{\det(X_t)}\right)^{\nu} Q_x^{m+\nu}|_{\mathscr{F}_t}.$$

**Proposition 7.** *For all $t > 0$ and $0 < \nu < 1$,*

$$Q_x^{m-\nu}(T_0 > t) = \frac{\Gamma_m(m)}{\Gamma_m(m+\nu)} \det\left(\frac{x}{2t}\right)^\nu {}_1F_1\left(\nu, m+\nu, -\frac{x}{2t}\right).$$

**Proof.** From the absolute-continuity relation above, we deduce that

$$Q_x^{m-\nu}(T_0 > t) = Q_x^{m+\nu}\left(\left(\frac{\det(x)}{\det(X_t)}\right)^\nu\right).$$

On the other hand, using the expression for the semi-group,

$$Q_x^\delta(\det(X_t)^s) = (2t)^{ms} \frac{\Gamma_m(s+\delta)}{\Gamma_m(\delta)} {}_1F_1\left(-s; \delta; -\frac{x}{2t}\right)$$

$$= (2t)^{ms} \frac{\Gamma_m(s+\delta)}{\Gamma_m(\delta)} \exp\left(-\operatorname{tr}\left(\frac{x}{2t}\right)\right) {}_1F_1\left(\delta+s; \delta; \frac{x}{2t}\right)$$

by the Kummer relation (cf. Theorem 7.4.3 in Muirhead [20]). Taking $s = -\nu$, we are done. $\square$

## 6. Generalized Hartman–Watson law

Henceforth, we will write $\mathscr{F}$ to denote one-dimensional hypergeometric functions. We define the *generalized Hartman–Watson* law as the law of

$$\int_0^t \operatorname{tr}(X_s^{-1}) \, ds \quad \text{under } Q_x^m(\cdot | X_t = y).$$

Its Laplace transform is given by

$$Q_x^m\left(\exp\left(\frac{-\nu^2}{2} \int_0^t \operatorname{tr}(X_s^{-1}) \, ds\right) \Big| X_t = y\right) = \frac{\Gamma_m(m)}{\Gamma_m(m+\nu)} \det(z)^{\nu/2} \frac{{}_0F_1(m+\nu, z)}{{}_0F_1(m, z)}, \quad (10)$$

$z = xy/4t^2$. Recall that for $m = 1$, this is the well-known Hartman–Watson law whose density was computed by Yor [23]. Here, we will investigate the case $m = 2$. The Gross–Richards formula is given, for $p = 0$ and $q = 1$, by

$${}_0F_1(m+\nu, z) = \frac{\det(z_i^{m-j} {}_0\mathscr{F}_1(m+\nu-j+1, z_i))}{V(z)},$$

where $(z_i)$ denote the eigenvalues of $z$ and $V(z) = \prod_{i<j}(z_i - z_j)$ is the Vandermonde determinant. Noting that $\Gamma_m(m+\nu) = \prod_{j=1}^m \Gamma(m+\nu-j+1)$, then

$$(10) = \frac{\det(z_i^{(m-j)/2} I_{m+\nu-j}(2\sqrt{z_i}))}{\det(z_i^{(m-j)/2} I_{m-j}(2\sqrt{z_i}))}.$$

Without loss of generality, we will take $t = 1$.

**Proposition 8.** *For $m = 2$, let $\lambda_1 > \lambda_2$ be the eigenvalues of $\sqrt{xy}$. Then the density of the generalized Hartman–Watson law is given by*

$$f(v) = \frac{\sqrt{\lambda_1 \lambda_2} v}{p\pi\sqrt{2\pi v^3}}$$

$$\times \frac{\int_0^1 \int_0^\infty z \sinh(p\sqrt{1-z^2}) e^{-2\sqrt{\lambda_1 \lambda_2} z \cosh y} e^{-2(y^2 - \pi^2)/v}(\sinh y) \sin(4\pi y/v) \, dz \, dy}{\int_0^1 \int_0^1 \frac{u \cosh(pu\sqrt{1-x^2})}{\sqrt{1-x^2}} I_0(2\sqrt{\lambda_1 \lambda_2} ux) \, du \, dx},$$

*for $v > 0$, where $p = \lambda_1 - \lambda_2$. Furthermore, if $\lambda_1 = \lambda_2 := \lambda$, then*

$$f(v) = \frac{4\lambda v e^{2\pi^2/v}}{\pi^2 \sqrt{2\pi v^3}} \frac{\int_0^\infty g(y) e^{-2y^2/v}(\sinh y) \sin(4\pi y/v) \, dy}{{}_1\mathscr{F}_2(\frac{1}{2}; 1; 2; \lambda^2)},$$

*where*

$$g(y) = \frac{1}{3} + \frac{\pi}{2} \frac{I_2(2\lambda \cosh y) + \mathbf{L}_2(2\lambda \cosh y)}{2\lambda \cosh y},$$

*and $\mathbf{L}_2$ is the Struve function Gradshteyn and Ryzhik [9].*

**Proof.** For $m = 2$, (10) becomes:

$$(10) = \frac{\lambda_1 I_{\nu+1}(\lambda_1) I_\nu(\lambda_2) - \lambda_2 I_{\nu+1}(\lambda_2) I_\nu(\lambda_1)}{\lambda_1 I_1(\lambda_1) I_0(\lambda_2) - \lambda_2 I_1(\lambda_2) I_0(\lambda_1)},$$

so, using the integral representations below (Brychkov *et al.* [3], page 46),

$$x(aI_{\nu+1}(ax) I_\nu(bx) - bI_{\nu+1}(bx) I_\nu(ax)) = (a^2 - b^2) \int_0^x u I_\nu(au) I_\nu(bu) \, du$$

with $x = 1, a = \lambda_1, b = \lambda_2$ and (Gradshteyn and Ryzhik [9], page 734)

$$\frac{\pi}{2} I_\nu\left(\frac{a}{2}(\sqrt{b^2 + c^2} + b)\right) I_\nu\left(\frac{a}{2}(\sqrt{b^2 + c^2} - b)\right) = \int_0^a \frac{\cosh(b\sqrt{a^2 - x^2})}{\sqrt{a^2 - x^2}} I_{2\nu}(cx) \, dx,$$

where $a > 0, \Re(\nu) > -1$. With $a = 1, b = (\lambda_1 - \lambda_2)u := pu$ and $c = 2\sqrt{\lambda_1 \lambda_2} u$, the numerator of (10) is then equal to:

$$\frac{2}{\pi}(\lambda_1^2 - \lambda_2^2) \int_0^1 \int_0^1 \frac{u \cosh(pu\sqrt{1-x^2})}{\sqrt{1-x^2}} I_{2\nu}(2\sqrt{\lambda_1 \lambda_2} ux) \, du \, dx.$$

Taking $\nu = 0$, the denominator is then equal to:

$$\frac{2}{\pi}(\lambda_1^2 - \lambda_2^2) \int_0^1 \int_0^1 \frac{u \cosh(pu\sqrt{1-x^2})}{\sqrt{1-x^2}} I_0(2\sqrt{\lambda_1 \lambda_2} ux) \, du \, dx.$$

Thus, (10) becomes

$$(10) = \frac{\int_0^1 \int_0^1 u\cosh(pu\sqrt{1-x^2})/\sqrt{1-x^2} I_{2\nu}(2\sqrt{\lambda_1\lambda_2}ux)\,du\,dx}{\int_0^1 \int_0^1 u\cosh(pu\sqrt{1-x^2})/\sqrt{1-x^2} I_0(2\sqrt{\lambda_1\lambda_2}ux)\,du\,dx}.$$

Now we only have to use the integral representation of $I_{2\nu}$ (Yor [23]):

$$I_{2\nu}(2\sqrt{\lambda_1\lambda_2}ux) = \frac{1}{2i\pi}\int_C e^{2\sqrt{\lambda_1\lambda_2}ux\cosh\omega}e^{-2\nu\omega}\,d\omega$$

$$= \frac{1}{2i\pi}\int_C e^{2\sqrt{\lambda_1\lambda_2}ux\cosh\omega}\int_0^\infty \frac{2\omega e^{-v\nu^2/2}}{(2\pi v^3)^{1/2}}e^{-2\omega^2/v}\,dv\,d\omega$$

where $C$ is the contour indicated in Yor [23]; hence, the density function is given by

$$f(v) = \frac{1}{i\pi\sqrt{2\pi v^3}}\frac{\int_0^1\int_0^1\int_C u\omega\cosh(pu\sqrt{1-x^2})/\sqrt{1-x^2}e^{2\sqrt{\lambda_1\lambda_2}ux\cosh\omega}e^{-2\omega^2/v}\,du\,dx\,d\omega}{\int_0^1\int_0^1 u\cosh(pu\sqrt{1-x^2})/\sqrt{1-x^2}I_0(2\sqrt{\lambda_1\lambda_2}ux)\,du\,dx}\mathbf{1}_{\{v>0\}}.$$

We can simplify this expression by integrating over $C$ to see that the numerator is equal to (Yor [23])

$$\frac{\sqrt{\lambda_1\lambda_2}v}{\pi\sqrt{2\pi v^3}}\int_0^1\int_0^1\int_0^\infty u^2 x\frac{\cosh(pu\sqrt{1-x^2})}{\sqrt{1-x^2}}e^{-2\sqrt{\lambda_1\lambda_2}ux\cosh y}$$

$$\times e^{-(-2(y^2-\pi^2))/v}(\sinh y)\sin\left(\frac{4\pi y}{v}\right)du\,dx\,dy.$$

Setting $z = ux$, the numerator is written

$$\frac{\sqrt{\lambda_1\lambda_2}v}{\pi\sqrt{2\pi v^3}}\int_0^1\int_0^u\int_0^\infty z\frac{u\cosh(p\sqrt{u^2-z^2})}{\sqrt{u^2-z^2}}e^{-2\sqrt{\lambda_1\lambda_2}z\cosh y}$$

$$\times e^{-2(y^2-\pi^2)/v}(\sinh y)\sin\left(\frac{4\pi y}{v}\right)du\,dz\,dy,$$

which we can integrate with respect to $u$ to obtain

$$\frac{\sqrt{\lambda_1\lambda_2}v}{p\pi\sqrt{2\pi v^3}}\int_0^1\int_0^\infty z\sinh(p\sqrt{1-z^2})e^{-2\sqrt{\lambda_1\lambda_2}z\cosh y}e^{-2(y^2-\pi^2)/v}(\sinh y)\sin\left(\frac{4\pi y}{v}\right)dz\,dy.$$

We now prove the second part. In this case $p = 0$ and we have to evaluate

$$\frac{\lambda v e^{2\pi^2/v}}{\pi\sqrt{2\pi v^3}}\frac{\int_0^1\int_0^1\int_0^\infty u^2 x/\sqrt{1-x^2}e^{-2\lambda ux\cosh y}e^{-2y^2/v}(\sinh y)\sin(4\pi y/v)\,du\,dx\,dy}{\int_0^1\int_0^1 uI_0(2\lambda ux)/\sqrt{1-x^2}\,du\,dx}.$$

Setting $z = ux$, the numerator reads

$$\frac{\lambda v e^{2\pi^2/v}}{\pi\sqrt{2\pi v^3}} \int_0^1 \int_0^\infty z\sqrt{1-z^2} e^{-2\lambda z \cosh y} e^{-2y^2/v} (\sinh y) \sin\left(\frac{4\pi y}{v}\right) dz\, dy.$$

Integration with respect to $z$ yields (Gradshteyn and Ryzhik [9], page 369):

$$\frac{\lambda v e^{2\pi^2/v}}{\pi\sqrt{2\pi v^3}} \int_0^\infty g(y) e^{-2y^2/v} (\sinh y) \sin\left(\frac{4\pi y}{v}\right) dy.$$

For the denominator, we use the fact that

$$\frac{d}{dz}(zI_1(z)) = zI_0(z),$$

which yields:

$$\int_0^1 \int_0^1 \frac{uI_0(2\lambda ux)}{\sqrt{1-x^2}} du\, dx = \int_0^1 \frac{I_1(2\lambda x)}{2\lambda x \sqrt{1-x^2}} dx.$$

Then, the formula

$$\int_0^a x^{\alpha-1}(a^2-x^2)^{\beta-1} I_\nu(cx) dx = 2^{-\nu-1} a^{2\beta+\alpha+\nu-2} c^\nu \frac{\Gamma(\beta)\Gamma((\alpha+\nu)/2)}{\Gamma(\beta+(\alpha+\nu)/2)\Gamma(\nu+1)}$$
$$\times {}_1\mathscr{F}_2\left(\frac{\alpha+\nu}{2}; \beta+\frac{\alpha+\nu}{2}; \nu+1; \frac{a^2c^2}{4}\right) \quad (11)$$

taken with $\alpha = 0, a = 1, \beta = 1/2, c = 2\lambda, \nu = 1$ gives

$$\int_0^1 \frac{I_1(2\lambda x)}{2\lambda x\sqrt{1-x^2}} dx = \frac{\pi}{4} {}_1\mathscr{F}_2\left(\frac{1}{2}; 1; 2; \lambda^2\right).$$

We can proceed differently. Let $\lambda_1 = \lambda_2 + h$. Then

$$(10) = \frac{((\lambda_2+h)I_{\nu+1}(\lambda_2+h)I_\nu(\lambda_2) - \lambda_2 I_{\nu+1}(\lambda_2) I_\nu(\lambda_2+h))/h}{((\lambda_2+h)I_1(\lambda_2+h)I_0(\lambda_2) - \lambda_2 I_1(\lambda_2) I_0(\lambda_2+h))/h}.$$

Next, we let $h \to 0$. As usual, we first compute the numerator and then take $\nu = 0$. To do this, we shall evaluate

$$A = \lim_{h\to 0} \frac{(\lambda_2+h)I_{\nu+1}(\lambda_2+h) - \lambda_2 I_{\nu+1}(\lambda_2)}{h},$$
$$B = \lim_{h\to 0} \frac{I_\nu(\lambda_2+h) - I_\nu(\lambda_2)}{h},$$

which are equal respectively to $\frac{d}{dx}(xI_{\nu+1}(x))$ and $\frac{d}{dx}(I_\nu(x))$ taken for $x = \lambda = \lambda_1 = \lambda_2$. Using the differentiation formula $\frac{d}{dx}(x^\nu I_\nu(x)) = x^\nu I_{\nu-1}(x)$ (Lebedev [18], page 110), we

obtain

$$\frac{\mathrm{d}}{\mathrm{d}x}(xI_{\nu+1}(x)) = -\nu I_{\nu+1}(x) + xI_\nu(x), \qquad \frac{\mathrm{d}}{\mathrm{d}x}(I_\nu(x)) = -\frac{\nu}{x}I_\nu(x) + I_{\nu-1}(x),$$

thus

$$N = I_\nu(\lambda)(-\nu I_{\nu+1}(\lambda) + \lambda I_\nu(\lambda)) - \lambda I_{\nu+1}(\lambda)\left(-\frac{\nu}{\lambda}I_\nu(\lambda) + I_{\nu-1}(\lambda)\right)$$

$$= \lambda(I_\nu^2(\lambda) - I_{\nu+1}(\lambda)I_{\nu-1}(\lambda)),$$

$$A = \frac{I_\nu^2(\lambda) - I_{\nu+1}(\lambda)I_{\nu-1}(\lambda)}{I_0^2(\lambda) - I_1(\lambda)I_{-1}(\lambda)}.$$

Using the integral representation (Gradshteyn and Ryzhik [9], page 757),

$$I_\mu(z)I_\nu(z) = \frac{2}{\pi}\int_0^{\pi/2} \cos((\mu-\nu)\theta)I_{\mu+\nu}(2z\cos\theta)\,\mathrm{d}\theta, \qquad \Re(\mu+\nu) > -1$$

the numerator is written as

$$N = \frac{2}{\pi}\int_0^{\pi/2}(1-\cos 2\theta)I_{2\nu}(2\lambda\cos\theta)\,\mathrm{d}\theta$$

$$= \frac{4}{\pi}\int_0^{\pi/2}(\sin^2\theta)I_{2\nu}(2\lambda\cos\theta)\,\mathrm{d}\theta$$

$$= \frac{4}{\pi}\int_0^1 \sqrt{1-r^2}I_{2\nu}(2\lambda r)\,\mathrm{d}r.$$

Thus, using (11), the denominator is equal to

$$D = \frac{4}{\pi}\int_0^1 \sqrt{1-r^2}I_0(2\lambda r)\,\mathrm{d}r = \frac{\pi}{4}{}_1\mathcal{F}_2\left(\frac{1}{2};1;2;\lambda^2\right).$$

Finally, the integral representation of $I_\nu$ gives

$$f(u) = \frac{\lambda u \mathrm{e}^{2\pi^2/u}}{\pi\sqrt{2\pi u^3}}\frac{\int_0^\infty \mathrm{e}^{-2y^2/u}\sinh(y)\sin(4\pi y/u)\int_0^1 r\sqrt{1-r^2}\mathrm{e}^{-2\lambda r\cosh y}\,\mathrm{d}r\,\mathrm{d}u}{\int_0^1 \sqrt{1-r^2}I_0(2\lambda r)\,\mathrm{d}r}$$

$$= \frac{\lambda u \mathrm{e}^{2\pi^2/u}}{\pi\sqrt{2\pi u^3}}\frac{\int_0^\infty g(y)\mathrm{e}^{-2y^2/u}\sinh(y)\sin(4\pi y/u)\,\mathrm{d}u}{\int_0^1 \sqrt{1-r^2}I_0(2\lambda r)\,\mathrm{d}r}. \qquad \square$$

## 7. The law of $T_0$

Recall that, for $0 < \nu < 1$,

$$Q_x^{m-\nu}(T_0 > t) = \frac{\Gamma_m(m)}{\Gamma_m(m+\nu)} \det\left(\frac{x}{2t}\right)^\nu {}_1F_1\left(\nu, m+\nu, -\frac{x}{2t}\right).$$

**Proposition 9.** *Let $m = 2$ and $\lambda_1 > \lambda_2$ be the eigenvalues of $x$. The density of $S_0 := 1/(2T_0)$ under $Q_x^{m-\nu}$ is given by*

$$f(u) = \frac{(\lambda_1 \lambda_2)^\nu u^{2\nu-2} e^{-(\lambda_1+\lambda_2)u}}{\Gamma(\nu+1)\Gamma(\nu)} \frac{{}_1\mathscr{F}_1(2, \nu+1, \lambda_1 u) - {}_1\mathscr{F}_1(2, \nu+1, \lambda_2 u)}{(\lambda_1 - \lambda_2)}.$$

**Corollary 7.** *If $\lambda_1 = \lambda_2 := \lambda$, the density is written*

$$f(u) = \frac{2\lambda^{2\nu} u^{2\nu-1} e^{-\lambda u}}{\Gamma(\nu+2)\Gamma(\nu)} {}_1\mathscr{F}_1(\nu-1, \nu+2, -\lambda u).$$

**Proof.** Recall first that when $m = 1$, $S_0 \stackrel{\mathcal{L}}{=} \gamma_\nu/x$, where $\gamma_\nu$ is a gamma variable with density $r^{\nu-1} e^{-r} dr$. With the help of the Gross–Richards formula, it follows that, for $m = 2$,

$$Q_x^{m-\nu}(S_0 \le u) = \frac{(\lambda_1 \lambda_2)^\nu}{(\lambda_1 - \lambda_2)\Gamma_2(\nu+2)} u^{2\nu}(\lambda_1 \, {}_1\mathscr{F}_1(\nu, \nu+2, -\lambda_1 u) \, {}_1\mathscr{F}_1(\nu-1, \nu+1, -\lambda_2 u)$$
$$- \lambda_2 \, {}_1\mathscr{F}_1(\nu, \nu+2, -\lambda_2 u) \, {}_1\mathscr{F}_1(\nu-1, \nu+1, -\lambda_1 u)),$$

where $S_0 := 1/(2T_0)$. This is a $C^\infty$ function in $u$. Hence, we will compute its derivative to obtain the density. Recall that

$$\frac{d}{dz} {}_1\mathscr{F}_1(a, b, z) = \frac{a}{b} {}_1\mathscr{F}_1(a+1, b+1, z),$$

thus

$$f(u) = \frac{d}{du} Q_x^{m-\nu}(S_0 \le u) = K(\nu, \lambda_1, \lambda_2) u^{2\nu-1}(A - B),$$

where

$$K(\nu, \lambda_1, \lambda_2) = \frac{(\lambda_1 \lambda_2)^\nu}{\Gamma_2(\nu+2)(\lambda_1 - \lambda_2)},$$

$$A = 2\nu(\lambda_1 \, {}_1\mathscr{F}_1(\nu, \nu+2, -\lambda_1 u) \, {}_1\mathscr{F}_1(\nu-1, \nu+1, -\lambda_2 u)$$
$$- \lambda_2 \, {}_1\mathscr{F}_1(\nu, \nu+2, -\lambda_2 u) \, {}_1\mathscr{F}_1(\nu-1, \nu+1, -\lambda_1 u)),$$

$$B = \frac{\nu}{\nu+2}(\lambda_1^2 u \, _1\mathscr{F}_1(\nu+1,\nu+3,-\lambda_1 u) \, _1\mathscr{F}_1(\nu-1,\nu+1,-\lambda_2 u)$$
$$- \lambda_2^2 u \, _1\mathscr{F}_1(\nu+1,\nu+3,-\lambda_2 u) \, _1\mathscr{F}_1(\nu-1,\nu+1,-\lambda_1 u)).$$

Then, we use the contiguous relation

$$b \, _1\mathscr{F}_1(a,b,z) - b \, _1\mathscr{F}_1(a-1,b,z) = z \, _1\mathscr{F}_1(a,b+1,z)$$

to see that

$$\lambda_1 u \, _1\mathscr{F}_1(\nu+1,\nu+3,-\lambda_1 u) = (\nu+2)(_1\mathscr{F}_1(\nu,\nu+2,-\lambda_1 u) - \, _1\mathscr{F}_1(\nu+1,\nu+2,-\lambda_1 u)),$$
$$\lambda_2 u \, _1\mathscr{F}_1(\nu+1,\nu+3,-\lambda_2 u) = (\nu+2)(_1\mathscr{F}_1(\nu,\nu+2,-\lambda_2 u) - \, _1\mathscr{F}_1(\nu+1,\nu+2,-\lambda_2 u))$$

imply that

$$f(u) = K_1(\nu,\lambda_1,\lambda_2) u^{2\nu-1}(C+D-E-F)$$

where

$$K_1(\nu,\lambda_1,\lambda_2) = \frac{\nu(\lambda_1\lambda_2)^\nu}{\Gamma_2(\nu+2)(\lambda_1-\lambda_2)},$$
$$C = \lambda_1 \, _1\mathscr{F}_1(\nu,\nu+2,-\lambda_1 u) \, _1\mathscr{F}_1(\nu-1,\nu+1,-\lambda_2 u),$$
$$D = \lambda_1 \, _1\mathscr{F}_1(\nu+1,\nu+2,-\lambda_1 u) \, _1\mathscr{F}_1(\nu-1,\nu+1,-\lambda_2 u),$$
$$E = \lambda_2 \, _1\mathscr{F}_1(\nu,\nu+2,-\lambda_2 u) \, _1\mathscr{F}_1(\nu-1,\nu+1,-\lambda_1 u),$$
$$F = \lambda_2 \, _1\mathscr{F}_1(\nu+1,\nu+2,-\lambda_2 u) \, _1\mathscr{F}_1(\nu-1,\nu+1,-\lambda_1 u).$$

Applying the above contiguous relation again yields:

$$\lambda_1 u \, _1\mathscr{F}_1(\nu+1,\nu+2,-\lambda_1 u) = (\nu+1)(_1\mathscr{F}_1(\nu,\nu+1,-\lambda_1 u) - \, _1\mathscr{F}_1(\nu+1,\nu+1,-\lambda_1 u)),$$
$$\lambda_2 u \, _1\mathscr{F}_1(\nu+1,\nu+2,-\lambda_2 u) = (\nu+1)(_1\mathscr{F}_1(\nu,\nu+1,-\lambda_2 u) - \, _1\mathscr{F}_1(\nu+1,\nu+1,-\lambda_2 u)),$$
$$\lambda_2 u \, _1\mathscr{F}_1(\nu,\nu+2,-\lambda_2 u) = (\nu+1)(_1\mathscr{F}_1(\nu-1,\nu+1,-\lambda_2 u) - \, _1\mathscr{F}_1(\nu,\nu+1,-\lambda_2 u)),$$
$$\lambda_1 u \, _1\mathscr{F}_1(\nu,\nu+2,-\lambda_1 u) = (\nu+1)(_1\mathscr{F}_1(\nu-1,\nu+1,-\lambda_1 u) - \, _1\mathscr{F}_1(\nu,\nu+1,-\lambda_1 u)).$$

Substituting in the expression for $f$, we obtain

$$f(u) = K_2(\nu,\lambda_1,\lambda_2) u^{2\nu-2}(G-H),$$

where

$$K_2(\nu,\lambda_1,\lambda_2) = \frac{\nu(\nu+1)(\lambda_1\lambda_2)^\nu}{\Gamma_2(\nu+2)(\lambda_1-\lambda_2)},$$
$$G = \, _1\mathscr{F}_1(\nu+1,\nu+1,-\lambda_2 u) \, _1\mathscr{F}_1(\nu-1,\nu+1,-\lambda_1 u),$$
$$H = \, _1\mathscr{F}_1(\nu+1,\nu+1,-\lambda_1 u) \, _1\mathscr{F}_1(\nu-1,\nu+1,-\lambda_2 u).$$

Eventually, writing

$$\Gamma_2(\nu+2) = \Gamma(\nu+2)\Gamma(\nu+1) = \nu(\nu+1)\Gamma(\nu+1)\Gamma(\nu),$$
$$_1\mathscr{F}_1(a,a,z) = e^{-z}\,_1\mathscr{F}_1(a,b,-z) = e^{-z}\,_1\mathscr{F}_1(b-a,b,z),$$

we obtain

$$f(u) = \frac{(\lambda_1\lambda_2)^\nu u^{2\nu-2} e^{-(\lambda_1+\lambda_2)u}}{\Gamma(\nu+1)\Gamma(\nu)} \frac{{}_1\mathscr{F}_1(2,\nu+1,\lambda_1 u) - {}_1\mathscr{F}_1(2,\nu+1,\lambda_2 u)}{\lambda_1 - \lambda_2}.$$

The case $\lambda_1 = \lambda_2$ is treated in the same way as before (for the Hartman–Watson law). In fact, writing $\lambda_1 = \lambda_2 + h$ and letting $h \to 0$, we see that the density is given by

$$\begin{aligned} f(u) &= \frac{\lambda^{2\nu} u^{2\nu-2} e^{-2\lambda u}}{\Gamma(\nu+1)\Gamma(\nu)} \frac{d}{d\lambda}\,_1\mathscr{F}_1(2,\nu+1,\lambda u) \\ &= \frac{2\lambda^{2\nu} u^{2\nu-1} e^{-2\lambda u}}{\Gamma(\nu+2)\Gamma(\nu)}\,_1\mathscr{F}_1(3,\nu+2,\lambda u) \\ &= \frac{2\lambda^{2\nu} u^{2\nu-1} e^{-\lambda u}}{\Gamma(\nu+2)\Gamma(\nu)}\,_1\mathscr{F}_1(\nu-1,\nu+2,-\lambda u). \end{aligned} \qquad \square$$

## 8. Conclusion

The Gross–Richards formula has been the main ingredient in this paper, since it enables us more explicitly to express the special functions with matrix argument. The case $m=3$ can be treated in the same way, but computation becomes too complicated. So, if we want to deal with the general case, it will be convenient to find a more explicit formula. Indeed, Schur functions can be expressed as polynomials in the elementary symmetric functions $e_r$ or as polynomials in the completely symmetric functions $h_r$. More precisely, we have

$$\begin{aligned} s_\lambda &= \det(e_{\lambda_i - i + j}), & 1 \leq i,j \leq n, \\ s_\lambda &= \det(h_{\lambda'_i - i + j}), & 1 \leq i,j \leq n, \end{aligned}$$

where $\lambda$ is a partition of length $\leq n$, and $\lambda'$ is the conjugate of $\lambda$ Macdonald [19]. So using these two identities, can we improve our results?

## Appendix: Special functions

### A.1. The hypergeometric series

The multivariate hypergeometric functions were studied by Muirhead [20] in the real symmetric case, Chikuze [4] for the complex Hermitian case, and Faraut and Korányi [8] in a more general setting. For Hermitian matrix argument, they are defined by

$$_pF_q((a_i)_{1\leq i\leq p}, (b_j)_{1\leq j\leq q}; X) = \sum_{k\geq 0} \sum_{\tau \perp k} \frac{(a_1)_\tau \cdots (a_p)_\tau}{(b_1)_\tau \cdots (b_q)_\tau} \frac{C_\tau(X)}{k!},$$

where $\tau = (k_1, \ldots, k_m)$ is a partition of weight $k$ and length $m$ such that $k_1 \geq \cdots \geq k_m$, $(a)_\tau$ is the generalized Pochammer symbol defined by

$$(a)_\tau = \prod_{i=1}^{m} \frac{\Gamma(a + k_i - i + 1)}{\Gamma(a - i + 1)}, \qquad \tau = (k_1, \ldots, k_m),$$

and $C_\tau(X)$ is the zonal polynomial of $X$ such that

$$(\operatorname{tr}(X))^k = \sum_{\tau \perp k} C_\tau(X).$$

Several normalizations for this polynomial exist in the litterature, but we do not consider them. This polynomial is symmetric, homogeneous, of degree $k$ in the eigenvalues of $X$, and is an eigenfunction of the differential operator

$$\Delta_X = \sum_{i=1}^{m} x_i^2 \frac{\partial^2}{\partial x_i^2} + 2 \sum_{i=1}^{m} \sum_{1 \leq k \neq i \leq m} \frac{x_i^2}{x_i - x_k} \frac{\partial}{\partial x_i}.$$

Furthermore, it is identified with the Schur function $s_\tau$ and $C_\tau(YX) = C_\tau(\sqrt{Y}X\sqrt{Y})$ for any Hermitian matrix $Y$. It is well known that, if $p = q + 1$, then the hypergeometric series is convergent for $0 \leq \|X\| < 1$ ($\|\cdot\|$ is the norm given by the spectral radius); if $p \leq q$, then it converges everywhere; otherwise it diverges.

### A.2. The modified Bessel function

The modified Bessel function Lebedev [18] with index $\nu \in \mathbb{R}$ is given by the series

$$I_\nu(z) = \sum_{k=0}^{\infty} \frac{1}{k!\Gamma(\nu + k + 1)} \left(\frac{z}{2}\right)^{2k+\nu}, \qquad z \in \mathbb{C}.$$

It can be represented through standard hypergeometric functions $_0F_1$ and $_1F_1$:

$$I_\nu(z) = \frac{1}{\Gamma(\nu + 1)} \left(\frac{z}{2}\right)^\nu {}_0F_1(\nu + 1; z^2).$$

## Acknowledgements

I would like to thank C. Donati-Martin (my doctoral thesis supervisor), J. Faraut and M. Yor for helpful comments and encouragement.